\documentclass[12pt]{amsart}
\usepackage{bbm}
\usepackage{mathrsfs}
\usepackage[english]{babel}
\usepackage{wrapfig}

\usepackage[top=3cm,bottom=2cm,left=3cm,right=3cm,marginparwidth=1.75cm]{geometry}

\usepackage{amsmath,amssymb}
\usepackage{graphicx}
\usepackage{cite}
\usepackage{enumitem} 
\usepackage[all,cmtip]{xy}
\usepackage{tikz-cd}
\usetikzlibrary{intersections}
\usepackage{relsize}
\usepackage{faktor}
\usepackage{dsfont}
\usepackage{xcolor, import}
\usepackage{bm}
\usepackage{hyperref}
\hypersetup{
 citebordercolor={green!40!black},
 linkbordercolor={red!50!black},
 urlbordercolor={blue!70!black}
}

\makeatletter
\renewcommand{\section}{\@startsection%
{section}% name
{1}% level
{0mm}% indent
{1.5\bigskipamount}% beforeskip
{0.5\bigskipamount}% afterskip
{\centering\normalsize\sc}}% style

\renewcommand{\subsection}{\@startsection%
{subsection}% name
{2}% level
{0mm}% indent
{0.5\bigskipamount}% beforeskip
{0.5mm}% afterskip
{\normalsize\sc}}% style

\renewcommand{\paragraph}{\@startsection%
{paragraph}% name
{4}% level
{0mm}% indent
{\bigskipamount}% beforeskip
{0pt}% afterskip
{\normalsize\bf}}% style

\expandafter\let\expandafter\oldproof\csname\string\proof\endcsname
\let\oldendproof\endproof
\renewenvironment{proof}[1][\proofname]{%
  \oldproof[\slshape #1]%
}{\oldendproof}
\def\provedboxcontents#1{$\square$}
\newtheoremstyle{thm}{6pt plus 1pt minus 1pt}{6pt plus 1pt minus 1pt}{\slshape}{}{\scshape}{.}{5pt plus 1pt minus 1pt}{}
\newtheoremstyle{def}{6pt plus 1pt minus 1pt}{6pt plus 1pt minus 1pt}{}{}{\scshape}{.}{5pt plus 1pt minus 1pt}{}
\newtheoremstyle{rmk}{6pt plus 1pt minus 1pt}{6pt plus 1pt minus 1pt}{}{}{\scshape}{.}{5pt plus 1pt minus 1pt}{}
\newtheoremstyle{claim}{6pt plus 1pt minus 1pt}{6pt plus 1pt minus 1pt}{}{}{\slshape}{.}{5pt plus 1pt minus 1pt}{}

\theoremstyle{thm}
\newtheorem{newstatement}{newstatement}
\newtheorem{lemma}[newstatement]{Lemma}
\newtheorem{theorem}[newstatement]{Theorem}
\newtheorem*{theorem*}{Theorem 2}

\newtheorem*{stiefel-thm}{Stiefel's Parallelizability Theorem}

\theoremstyle{def}
\newtheorem{definition}[newstatement]{Definition}

\theoremstyle{rmk}
\newtheorem{remark}[newstatement]{Remark}
\newtheorem{example}[newstatement]{Example}
\newtheorem*{example*}{Example}

\theoremstyle{claim}

\renewcommand{\epsilon}{\varepsilon}
\renewcommand{\phi}{\varphi}

\newcommand{\R}{\mathbb{R}}

\newcommand{\N}{\mathbb{N}}
\newcommand{\Z}{\mathbb{Z}}

\DeclareMathOperator{\cs}{\#}
\newcommand{\simtimes}{\mathbin{\widetilde{\smash{\times}}}}

\let\emph\textsl
\title{$\text{Pin}^{\pm}$-structures on non-oriented 4-manifolds via Lefschetz fibrations}

\author{Valentina Bais}
\address{Department of Mathematics, SISSA, Via Bonomea 265, 34136 Trieste, Italy.}
\email{vbais@sissa.it}

\date{\today}

\begin{document}

%\normalpage

\begin{abstract}
We study necessary and sufficient conditions for a 4-dimensional Lefschetz fibration over the 2-disk to admit a $\text{Pin}^{\pm}$-structure, extending the work of A. Stipsicz in the orientable setting. As a corollary, we get existence results of $\text{Pin}^{+}$ and $\text{Pin}^-$-structures on closed non-orientable 4-manifolds and on Lefschetz fibrations over the 2-sphere. In particular, we show via three explicit examples how to read-off $\text{Pin}^{\pm}$-structures from the Kirby diagram of a 4-manifold. We also provide a proof of the well-known fact that any closed 3-manifold $M$ admits a $\text{Pin}^-$-structure and we find a criterion to check whether or not it admits a $\text{Pin}^+$-structure in terms of a handlebody decomposition. We conclude the paper with a characterization of $\text{Pin}^+$-structures on vector bundles.
\end{abstract}


\maketitle

\section{Introduction}

$\text{Pin}^{\pm}$-structures can be thought of as the non-oriented analogue of $\text{Spin}$-structures. We refer the reader to \cite{[KirbyTaylor]} and \cite[Section 2]{[Stolz]} for a precise definition of such Lie groups and for the background on $\text{Pin}^{\pm}(n)$-structures on low dimensional manifolds. In this note, we find necessary and sufficient conditions for a (possibly non-orientable) Lefschetz fibration over the 2-disk to support a $\text{Pin}^{-}$ or a $\text{Pin}^+$-structure. In particular, such conditions are expressed in terms of the homology classes of the vanishing cycles of the Lefschetz fibration, with coefficients taken in $\Z_2$ and $\Z_4$ respectively. We refer the reader to Section \ref{background} for a quick recap on basic facts about non-orientable Lefschetz fibrations. Our main reference for this subject is \cite{[MillerOzbagci]}. 

In Section \ref{Pin 1} we prove the following statement, extending the work of A. Stipsicz on $\text{Spin}$-structures over Lefschetz fibrations in \cite[Theorem 1.1]{[Stipsicz]}.  

\begin{theorem}\label{theorem1}
    Let $X$ be a smooth 4-manifold and $f: X \to D^2$ a Lefschetz fibration with regular fiber $\Sigma$. There is no $\text{Pin}^-$-structure on $X$ if and only if there are $k+1$ vanishing cycles $c_0, c_1, \dots, c_k$ such that  $[c_0]=\sum_{i=1}^k [c_i] \in H_1(\Sigma; \Z_2)$ and $k+\sum_{1 \leq i < j \leq k} c_i \cdot c_j \equiv 0 \mod{2}$. 
\end{theorem}

As we will see in Section \ref{Pin 1}, Theorem \ref{theorem1} is a consequence of the fact that $\text{Pin}^-$- structures on a Lefschetz fibration over the 2-disk with regular fiber $\Sigma$ naturally correspond to maps
\[ q^-: H_1(\Sigma;\Z_2) \rightarrow \Z_4 \]
such that $q^-(x+y)=q^-(x)+q^-(y)+2x \cdot y$ for every $x, y \in H_1(\Sigma;\Z_2)$ and $q^-([c])=2$ on every vanishing cycle $c \subset \Sigma$, see Lemma \ref{lemma1}. Moreover, such correspondence is equivariant with respect to the action of $H^1(X;\Z_2)$.

Section \ref{Pin 2} is devoted to the study of $\text{Pin}^+$-structures on Lefschetz fibrations. By \cite{[D]}, a $\text{Pin}^+$-structure on the regular fiber $\Sigma$ corresponds to a map
\[q_0^+: H_1(\Sigma;\Z_4) \rightarrow \Z_2\]
with the property that $q_0^+(x+y)=q_0^+(x)+ q_0^+(y)+x \cdot y \in \Z_2$ for all $x,y \in H_1(\Sigma;\Z_4)$. We show that we can choose $q_0^+$ to be the restriction of a $\text{Pin}^+$-structure defined on the whole fibration if and only if the following holds.

\begin{theorem}\label{Main2}
The total space of the Lefschetz fibration $f:X \to D^2$ with vanishing cycles $c_1, \dots, c_n \subset \Sigma$ supports a $\text{Pin}^+$-structure if and only if $\Sigma$ supports a $\text{Pin}^+$-structure and \[\text{rank} (C)= \text{rank}(C \mid A)\]
    where $C$ is the $\Z_2$-reduction of the $n \times r$ matrix $(c_{ij})$ whose rows are given by the components of $[c_1], \dots, [c_n] \in H_1(\Sigma;\Z_4) \cong \Z^r$ with respect to a fixed basis $e_1, \dots, e_r$ of the free $\Z_4$-module $H_1(\Sigma;\Z_4)$ and $A$ is the column vector with entries $A_i=1+q_0^+([c_i]) \in \Z_2$ for $i=1, \dots, n$.
\end{theorem}

A key tool in the proof of the above result is Lemma \ref{lemma2}, in which we show that $\text{Pin}^+$-structures on $X$ correspond to maps
\[q^+: H_1(\Sigma;\Z_4) \rightarrow \Z_2\]
such that $q^+(x+y)=q^+(x)+q^+(y)+x\cdot y$ for every $x,y \in H_1(\Sigma;\Z_4)$ and $q^+([c])=1$ on every vanishing cycle $c \subset \Sigma$. The correspondence is also in this case equivariant with respect to the $H^1(X;\Z_2)$-action.

We remark that the study of $\text{Pin}^{+}$ and $\text{Pin}^-$-structures on non-orientable 4-manifolds is essentially reduced to the case of non-orientable Lefschetz fibrations over the 2-disk due to the following result, which is a straightforward consequence of \cite[Theorem 1.1]{[MillerOzbagci]}.

\begin{theorem}[Miller-Ozbagci \cite{[MillerOzbagci]}]\label{T}
Let $X$ be a closed non-orientable smooth 4-manifold. There is a decomposition \[X=L \cup_{\partial} H\] where $L$ is a non-orientable Lefschetz fibration over the 2-disk and $H$ is a non-orientable $1$-handlebody. 
\end{theorem}

In particular, the proof of \cite[Theorem 1.1]{[MillerOzbagci]} shows an explicit way of endowing a (possibly non-orientable) 2-handlebody with the structure of a Lefschetz fibration over the 2-disk, see also \cite{[Harer]} and \cite{[EF]}. Moreover, a 4-manifold $X=L \cup_{\partial} H$ as in Theorem \ref{T} admits a $\text{Pin}^{\pm}$-structure if and only if $L$ does. This latter conditions can be checked using Theorem \ref{theorem1} and Theorem \ref{Main2}. 

In Section \ref{examples} we provide explicit examples in which, starting from a handlebody decomposition of a non-orientable 4-manifold, we endow its 2-handlebody with the structure of a Lefschetz fibration over the 2-disk. We then apply our results in order to check whether such 4-manifolds support a $\text{Pin}^-$ or a $\text{Pin}^+$-structure. In the case they do, we find all the possible structures by looking at the $H^1(X,\Z_2)$-action on the associated quadratic enhancements on the regular fibers. The non-orientable handlebodies under consideration describe the three 4-manifolds $\R \mathbb{P}^4$, $S^2 \simtimes \R \mathbb{P}^2$ and $S^2 \times \R \mathbb{P}^2$ and can be found in \cite{[Akbulut1]}, \cite{[BaisTorres]}, \cite{[MillerNaylor]} and \cite{[Torres]}. In particular, that there are two non-orientable $S^2$-bundles over $\R \mathbb{P}^2$ and we denote by $S^2 \simtimes \R \mathbb{P}^2$ the non-trivial one, see \cite[Section 2.3]{[BaisTorres]}.

In analogy to the $\text{Spin}$ case, one can define the $\text{Pin}^{+}$ and $\text{Pin}^-$ cobordism groups in each dimension, see \cite[Section 1]{[KirbyTaylor]}. As sets, they consist of the equivalence classes of closed $n$-dimensional smooth manifolds with a fixed $\text{Pin}^{+}$ (resp. $\text{Pin}^-$)-structure, which are identified whenever they co-bound a $\text{Pin}^{+}$ (resp. $\text{Pin}^-$) $(n+1)$-manifold. The group operation is the one induced by the disjoint union of manifolds. We remark that $\text{Pin}^+$ and $\text{Pin}^-$ cobordisms groups are drastically different. In particular, the 4-dimensional $\text{Pin}^+$-cobordism group is 
\[ \Omega_4^{\text{Pin}^+} \cong \Z_{16}\]
(see \cite[Section 3]{[KirbyTaylor]}) and the $\eta$-invariant modulo $2 \Z$ of the twisted Dirac operator associated to a $\text{Pin}^+$-structure on a 4-manifold is a complete invariant for $\text{Pin}^+$-cobordism classes, see \cite{[Stolz]}. In particular, there are instances in which such invariant can detect exotic behaviors in the non-orientable 4-dimensional realm, in the sense that there are examples of non-orientable $\text{Pin}^+$ 4-manifolds which share the same homeomorphism type but define two different classes in $\Omega_4^{\text{Pin}^+}$ and hence cannot be diffeomorphic. To the best of our knowledge, following a suggestion in \cite{[eta]}, this was shown for the first time in \cite{[Stolz]}, where the $\eta$-invariant of $\text{Pin}^+$-structures is used to detect the exotic $\R \mathbb{P}^4$ constructed in \cite{[CappellShaneson]}. We remark that another exotic $\R \mathbb{P}^4$ is constructed in \cite{[FS]} and its $\eta$-invariant is computed in \cite{[O2]}. An analogous approach has then been used also in \cite{[BaisTorres]} and \cite{[Torres]}. On the other hand
\[\Omega_4^{\text{Pin}^-} =0\] 
(see again \cite[Section 3]{[KirbyTaylor]}) and hence it is not possible to find $\text{Pin}^-$ exotic 4-manifolds using this strategy. In particular, the study of $\text{Pin}^+$-structures on 4-manifolds could potentially lead to  better understanding of exotic behaviors in the non-orientable realm and this fact is one of the main motivations for this note. In particular, in Section \ref{Pin 2} we show that if $X =L \cup H$ is a closed 4-manifold with a decomposition as in the statement of Theorem \ref{T}, then the datum of a $\text{Pin}^+$-structure on $X$ is equivalent to the datum of a map \[q^+: H_1(\Sigma;\Z_4)\rightarrow \Z_2\] such that $q^+(x+y)=q^+(x)+q^+(y)+ x \cdot y$ for all $x,y \in H_1(\Sigma;\Z_4)$ taking value $1 \in \Z_2$ on all the vanishing cycles of $L$, where $\Sigma$ denotes a regular fiber. We leave the reader with the following question.

    \textbf{Question}: Is there a formula for the $\eta$-invariant modulo $2 \Z$ of a $\text{Pin}^+$-structure on $X$ in terms of the corresponding map $q^{+}:H_1(\Sigma;\Z_4) \rightarrow \Z_2$?

Section \ref{old} contains a short discussion on how to check the vanishing of $w_2$ and $w_1^2$ is terms of the embeddings of homologically essential surfaces inside the ambient 4-manifold. In Section \ref{sphere} we study Lefschetz fibrations over the 2-sphere, \ref{3-maifolds} is devoted to a brief discussion on the 3-dimensional case, while Section \ref{Milnor} contains a characterization of $\text{Pin}^+$-structures on vector bundle which resembles Milnor's characterization of $\text{Spin}$-structures, see \cite[Alternative definitions 2]{[Milnor]}.

All manifolds in this note are smooth. We will not assume our 4-manifolds to be orientable, unless otherwise stated.

\section*{Acknowledgments}
I am grateful to my advisors Rafael Torres and Daniele Zuddas for the constant support and for their comments on this manuscript. 
The author has been partially supported by GNSAGA – Istituto Nazionale di Alta Matematica ‘Francesco Severi’, Italy. 
\section{$\text{Pin}^{\pm}$-structures and embedded surfaces}\label{old}

 For every $n \in \N$, there are two central extensions 
\begin{equation}\label{cover}
    \text{Pin}^{\pm}(n) \rightarrow O(n)
\end{equation}
which are topologically $\text{Spin}(n) \cup \text{Spin}(n)$ but algebraically not, in the sense that they are both semi-direct products of the form $\text{Spin}(n) \rtimes_{\phi_{\pm}} \Z_2$ for different choice of homomorphism $\phi_{\pm}: \Z_2 \rightarrow \text{Aut}(\text{Spin(n)})$, see \cite[Proposition 6.2]{[Stolz]}. In the orientable setting, it is known that the existence of a $\text{Spin}$-structure on a real vector bundle $\xi$ on a manifold $X$ is equivalent to the vanishing of $w_2(\xi) \in H^2(X,\Z_2)$, where $w_i$ is the $i^{th}$ Stiefel-Whitney class of $\xi$. This equivalent to having an even intersection form, provided that  $H^2(M;\Z)$ has no 2-torsion \cite[Corollary 5.7.6]{[GompfStipsicz]}. This can be shown via the Wu formula, which states that for any $a \in H_2(M;\Z_2)$ one has \begin{equation}\label{Wu}\langle w_2(M),a\rangle = a^2 \text{ modulo} \ 2\end{equation} see \cite[Proposition 1.4.18]{[GompfStipsicz]}. In a similar way, there are cohomological obstructions to the existence of $\text{Pin}^+$ and $\text{Pin}^-$-structures on the tangent bundle of a manifold, namely the vanishing of $w_2(\xi)$ and of $(w_2+w_1^2)(\xi)$ respectively, see \cite[Lemma 1.3]{[KirbyTaylor]}. Moreover, both $\text{Pin}^{+}$ and $\text{Pin}^-$-structures (when they do exist) are torsors over $H^1(X,\Z_2)$, meaning that there is a free and transitive action of $H^1(X;\Z_2)$ over the sets of such geometric structures. However, the Wu formula does not hold when $M$ is non-orientable, but one can still interpret the vanishing of $w_2(M)$ and $(w_2+w_1^2)(M)$ in terms of properties of embedded surfaces. Indeed, for any 4-manifold $M$ and $a \in H_2(M;\Z_2)$ there is a (possibly non-orientable) smoothly embedded surface $\Sigma \subset M$ representing $a$ (see \cite[Remark 1.2.4]{[GompfStipsicz]}) and one can show that
\begin{equation}\label{+}
    \langle w_2(M),a\rangle = \chi(\Sigma)  +(w_1(\Sigma) \smile w_1(\nu(\Sigma)))([\Sigma])+[\Sigma]^2 \text{ modulo } 2.
\end{equation}

Since the tangent bundle of a 4-manifold supports a $\text{Pin}^+$-structure if and only if $w_2(M)$ vanishes, this is equivalent to check that the evaluation (\ref{+}) is trivial on any homology class $a \in H_2(M;\Z_2)$.

On the other hand, the obstruction for $\text{Pin}^-$-structures to exist is $(w_2+w_1^2)(M)$ \cite[Lemma 1.3]{[KirbyTaylor]} so, in order to understand when such structures exist, we need to compute $w_1^2(M)$ and one can easily show that
\begin{equation}\label{-}
    \langle w_1^2(M),[\Sigma] \rangle = w_1^2(\Sigma)+ w_1^2(\nu(\Sigma)).
\end{equation}

\begin{example}[$\R \mathbb{P}^4$] The 4-dimensional real projective space $\R \mathbb{P}^4$ is an example of non-orientable 4-manifold supporting two distinct $\text{Pin}^+$-structures but no $\text{Pin}^-$-structure. We have that $H_2(\R \mathbb{P}^4;\Z_2) \cong \Z_2$ is generated by the homology class of $\R \mathbb{P}^2$. Moreover, the tubular neighbourhood $\nu(\R \mathbb{P}^2) \subset \R \mathbb{P}^4$ is diffeomorphic to the twisted 2-disk bundle $D^2 \simtimes \R \mathbb{P}^2$ defined as the quotient of $D^2\times S^2$ via the involution
\begin{align*}
    D^2 \times S^2 \rightarrow D^2 \times S^2
    \\ (x,y) \mapsto (\rho_{\pi}(x), -y)
\end{align*}
where $\rho_{\pi}$ denotes the rotation of $S^2$ of $\pi$ radians about a fixed axes.
One can show that such a quotient is diffeomorphic to
\begin{equation}
    D^2 \simtimes \R \mathbb{P}^2 \cong (D^2 \times D^2) \cup_{\phi} (D^2 \times \text{Mb}) 
\end{equation}
    where $\text{Mb}$ denotes the Möbius strip and $\phi$ is the map
    \begin{equation*}
        \phi: D^2 \times S^1 \rightarrow D^2 \times S^1, \quad (x, \theta) \mapsto (\rho_{\theta}(x), \theta)
    \end{equation*}
    and $\rho_{\theta}: D^2 \rightarrow D^2$ denotes the rotation of the 2-disk of angle $\theta$ with respect to the origin, see \cite[Section 0]{[Akbulutbook]} and \cite[Example 7]{[BaisTorres]}.

    By (\ref{+}) we have that
    \begin{equation*}
        \langle w_2(\R \mathbb{P}^4), [\R \mathbb{P}^2] \rangle = \chi(\R \mathbb{P})^2+ (w_1(\R \mathbb{P}^2)\cup w_1(\nu (\R \mathbb{P}^2)))([\R \mathbb{P}^2])+[\R \mathbb{P}^2]^2=0 \in \Z_2
    \end{equation*}
    since $\chi(\R \mathbb{P}^2)=1=[\R \mathbb{P}^2]^2$ and $w_1(\nu (\R \mathbb{P}^2))=0$, while (\ref{-}) implies that
    \begin{equation*}
        \langle w_1(\R \mathbb{P}^4)^2, [\R \mathbb{P}^2] \rangle = w_1^2(\R \mathbb{P}^2) + w_1^2(\nu (\R \mathbb{P}^2))= 1+0=1\in \Z_2
    \end{equation*}
    where $w_1^2(\R \mathbb{P}^2)$ is computed as the self-intersection of a loop in $\R \mathbb{P}^2$ whose $\Z_2$-homology class generates $H_1(\R \mathbb{P}^2;\Z_2)$.
\end{example}

\section{Non-orientable Lefschetz fibrations}\label{background}

Since the conditions \ref{+} and \ref{-} are not immediate to check in the non-orientable setting and depend heavily on the embeddings of the homologically essential surfaces, in the following sections we develop another way to combinatorially understand when a 4-manifold is $\text{Pin}^+$ or $\text{Pin}^-$ by means of a specific kind of decomposition. To do this, we will need the notion of non-orientable Lefschetz fibration over the 2-disk.

We start by recalling the definition of Lefschetz fibration. Our main reference for this topic in the non-orientable realm is \cite{[MillerOzbagci]}.

\begin{definition}
Let $X$ be a compact, connected 4-manifold and let $B$
be a compact, connected surface, both with possibly non-empty boundary. A Lefschetz fibration is a smooth submersion 
\[f : X \rightarrow B \]
away from finitely many points in the interior of $B$ such that each fiber contains at most one critical point and $f$ is a fiber bundle with surface fiber on the complement of the critical values. Moreover, around each critical point, we require $f$ to conform to the local complex model \[(z_1, z_2) \rightarrow z_1z_2.\]

\end{definition}

\begin{remark}
    If $X$ is non-orientable and $B$ is oriented, each regular fiber is a non-orientable surface. In this case it makes no sense to ask for the orientation of the local model around a critical point to be compatible to the one of $X$.
\end{remark}
In the following, we will just consider the case in which $B=D^2$ and $\Sigma$ will be regular fiber. It is possible to show that every Lefschetz fibration
\[f: X \rightarrow D^2\]
is obtained from the product one
\[\Sigma \times D^2 \rightarrow D^2\]
by gluing 4-dimensional 2-handles along $\Sigma \times \partial D^2$ with framing $\pm 1$ with respect to the fiber framing. The attaching curves of such 2-handles are push-offs in distinct fibers of simple closed curves inside $\Sigma$, which are called the vanishing cycles of the fibration and are denoted in the following by $c_1, \dots, c_n \subset \Sigma$. Moreover, $\langle w_1(X),[c_i] \rangle =0$ for all $i=1, \dots, n$ and hence the tubular neighbourhood of these attaching regions is necessarily diffeomorphic to a product of $S^1$ with the unit interval. We will call \textit{two-sided} all curves with this property.

\section{$\text{Pin}^-$-structures on 4-manifolds via Lefschetz fibrations}\label{Pin 1}

It is well known that the datum of a $\text{Spin}$-structure over an orientable surface is equivalent to the one of a quadratic enhancement \[ s: H_1(\Sigma;\Z_2) \to \Z_2\] satisfying the condition \[s(x+y)=s(x)+s(y)+ x \cdot y\] for all $x,y \in H_1(\Sigma;\Z_2)$, where $\cdot$ denotes the $\Z_2$-intersection number between cycles, see \cite{[Johnson]} and \cite{[Stipsicz]}.  A geometric interpretation of this algebraic object can be given as follows. The 1-dimensional $\text{Spin}$-cobordism group is $\Omega_1^{\text{Spin}}\cong \Z_2$ \cite{[Milnor]} and every closed simple curve $\gamma \subset \Sigma$ in a $\text{Spin}$ surface inherits a $\text{Spin}$-structure. In particular, if $s$ is the quadratic enhancement associated to the fixed $\text{Spin}$-structure on $\Sigma$, then $s([\gamma])=0$ if and only if the induced $\text{Spin}$-structure on $\gamma$ is the one bounding the unique $\text{Spin}$-structure on the 2-disk.

It is possible to give a similar description for $\text{Pin}^-$ structures on (not necessarily orientable) surfaces, as shown by the following result.

\begin{theorem}[Kirby-Taylor, \cite{[KirbyTaylor]}]\label{pin}
    There is a 1:1 correspondence between $\text{Pin}^-$-structures on a surface $\Sigma$ and quadratic enhancements \[ q^{-}: H_1(\Sigma; \Z_2) \to \Z_4\] with the property that \[q^{-}(x+y)=q^{-}(x)+q^{-}(y)+2 \ x \cdot y\] for any $x,y \in H_1(\Sigma; \Z_2)$.
\end{theorem}
In particular, if $\gamma \subset \Sigma$ is a simple closed curve, then $q^{-}([\gamma])$ is even if and only if $\gamma$ is two-sided and this is the case whenever $\gamma$ is a vanishing cycle for a Lefschetz fibration. Moreover, recall that the 1-dimensional $\text{Pin}^-$-cobordism group is $\Omega_1^{\text{Pin}^-}\cong \Z_2$ \cite[Section 0]{[KirbyTaylor]} and every closed simple curve $\gamma \subset \Sigma$ in a $\text{Pin}^-$ surface inherits a $\text{Pin}^-$ structure. We have that $q^{-}([\gamma])=0$ if and only if $\gamma$ inherits the bounding $\text{Pin}^-$-structure.

In order to prove Theorem \ref{theorem1}, we will need the following Lemma.
\begin{lemma}\label{lemma1}
    Let $f: X \rightarrow D^2$ be a Lefschetz fibration with regular fiber $\Sigma$ and vanishing cycles $c_1, \dots, c_n \subset \Sigma$. There is a natural one to one correspondence between the set of $\text{Pin}^-$-structures on $X$ and the set of quadratic enhancements \[q^-:H_1(\Sigma;\Z_2) \rightarrow \Z_4\] such that $q^-(x+y)=q^-(x)+q^-(y)+2x\cdot y$ for all $x, y \in H_1(\Sigma;\Z_2)$ and $q^-([c_i])=2$ for all $i=1, \dots, n$. Moreover, such correspondence is equivariant with respect to the free and transitive action of $H^1(X;\Z_2)$.
\end{lemma}

\begin{proof}
    The existence of a $\text{Pin}^-$-structure on $X$ is equivalent to the existence of a $\text{Pin}^-$-structure on $\Sigma$ that extends to the 2-handles. This follows from the fact that any $\text{Pin}^-$-structure on $X$ induces by restriction a $\text{Pin}^-$ structure on a (trivial) tubular neighbourhood $\nu(\Sigma)\cong \Sigma \times D^2$ and all the $\text{Pin}^-$-structures on $\Sigma \times D^2$ are pull-backs of the ones on $\Sigma$ via the projection map $\Sigma \times D^2 \rightarrow \Sigma$. Since such 2-handles are attached to $\Sigma \times \partial D^2$ with relative odd framing, this corresponds to a $\text{Pin}^-$-structure on $\Sigma$ restricting to the non-bounding one on every vanishing cycle. The conclusion follows from Theorem \ref{pin}.

     The free and transitive action of $H^1(X;\Z_2)$ on the set of $\text{Pin}^-$-structures on $X$ can be seen as follows. In \cite[Section 3]{[KirbyTaylor]} it is shown that $H^1(\Sigma;\Z_2)$ acts on the set of $\text{Pin}^-$-structures of the surface $\Sigma$ by
    \begin{equation}\label{act}
        q^{-}_{\gamma}(x)=q^{-}(x)+2 \cdot \gamma(x)
    \end{equation}
    for all $\gamma \in H^1(\Sigma;\Z_2)$ and $x \in H_1(\Sigma;\Z_2)$. In particular, $\text{Pin}^-$-structures on $\Sigma$ equivariantly correspond to quadratic enhancements and the action (\ref{act}) is well defined on $H^1(X;\Z_2)\subset H^1(\Sigma;\Z_2)\cong H^1(X;\Z_2)\oplus K$, where $K \subset H^1(\Sigma;\Z_2)$ is the kernel of the map induced by the inclusion $\Sigma \subset X$.
\end{proof}
In light of this fact, it is possible to characterize the existence of $\text{Pin}^-$-structures on non-orientable Lefschetz fibrations over the 2-disk in terms of the $\Z_2$-homology classes of the vanishing cycles. 

\begin{proof}[Proof of Theorem \ref{theorem1}]
    By Lemma \ref{lemma1}, $X$ does not support any $\text{Pin}^-$-structure if and only if it is not possible to build a map
    \[q^-: H_1(\Sigma; \Z_2) \rightarrow \Z_4 \]
    such that $q^-(x+y)=q^-(x)+q^-(y)+2x\cdot y$ for all $x,y \in H_1(\Sigma;\Z_2)$ and such that $q^+([c_i])=2$ for all $i=1, \dots, n$. The remaining part of the proof is essentially the same as the proof of \cite[Theorem 1.1]{[Stipsicz]}. The details are left to the interested reader. 
\end{proof}

Note that, if we restrict to orientable surfaces, there is a natural map \[\Omega_1^{\text{Spin} }\to \Omega_1^{\text{Pin}^- }\] giving a group isomorphism, see \cite[Theorem 2.1]{[KirbyTaylor]}. At the level of the associated quadratic forms, the enhancement
\[s: H_1(\Sigma;\Z_2)\to \Z_2\]
corresponds to
\[q^{-}: H_1(\Sigma;\Z_2)\to \Z_4\]
given by\[q^{-}(x) = 2 \cdot s(x) \]
for every $x \in H_1(\Sigma;\Z_2)$, where $2 \cdot$ denotes the inclusion $\Z_2 \subset \Z_4$. In particular, the condition on the vanishing cycles we found in Theorem \ref{theorem1} coincides with the one in \cite[Theorem 1.1]{[Stipsicz]} and this is due to the fact that, when restricting to orientable Lesfschetz fibrations over $D^2$, being $\text{Spin}$ coincides with being $\text{Pin}^-$.

\section{$\text{Pin}^+$-structures on Lefschetz fibrations over the 2-disk}\label{Pin 2}

As a consequence of \cite[Theorem A]{[D]}, there is a canonical affine bijective correspondence between $\text{Pin}^+$-structures on a surface $\Sigma$ and the set of maps 
\begin{equation}
    q^{+}: H_1(\Sigma;\Z_4) \to \Z_2
\end{equation}
with the property that 
\begin{equation}\label{Rule}
    q^{+}(x+y)=q^{+}(x)+q^{+}(y)+x \cdot y
\end{equation}
  for every $x, y \in H_1(\Sigma,\Z_4)$. 

Recall that the 1-dimensional $\text{Pin}^+$-cobordism group is
\[\Omega_1^{\text{Pin}^+}\cong \{0\}\]
see \cite[Section 3]{[KirbyTaylor]}.
However, $S^1$ has two distinct $\text{Pin}^+$-structures. These correspond to two different trivializations of the direct sum $TS^1 \oplus \epsilon$ of the tangent bundle of $S^1$ with a trivial line bundle, where the first one is given by the restriction to $S^1$ of a trivialization of the tangent bundle $T D^2$ of the 2-disk and the second one comes from the Lie group framing of $S^1$. In particular, one can show that such structures have the property that they respectively bound the 2-disk and the Möbius strip, see \cite[Section 1]{[KirbyTaylor]}. A close look at the construction of $q^{+}$ starting from a $\text{Pin}^+$-structure on $\Sigma$ shows that $q^{+}([\gamma])=0$ if the induced $\text{Pin}^+$-structure on $\gamma$ is the one bounding a 2-disk, while $q^{+}([\gamma])=1$ if it is the one bounding the Möbius strip, see \cite{[KirbyTaylor]}. Note that, given any two maps $q^{+}_1, q^{+}_2$ corresponding to $\text{Pin}^+$-structures on $\Sigma$, their difference can be regarded as an element \[q^{+}_1 - q^{+}_2 \in\text{Hom}(H_1(\Sigma;\Z_4),\Z_2).\]

The proof of Theorem \ref{Main2} is based on the following Lemma.
\begin{lemma}\label{lemma2}
    Let $f: X \rightarrow D^2$ be a Lefschetz fibration with regular fiber $\Sigma$ and vanishing cycles $c_1, \dots, c_n \subset \Sigma$. There is a natural bijection between the set of $\text{Pin}^+$-structures on $X$ and the set of quadratic enhancements \[q^+:H_1(\Sigma;\Z_4) \rightarrow \Z_2\] such that
    $q^+(x+y)=q^+(x)+q^+(y)+x \cdot y$ for all $x, y \in H_1(\Sigma; \Z_4)$ and $q^+([c_i])=1$ for all $i=1, \dots, n$. Moreover, such correspondence is equivariant with respect to the free and transitive action of $H^1(X;\Z_2)$.
\end{lemma}

\begin{proof}
    As is the $\text{Pin}^-$ case, the existence of a $\text{Pin}^+$-structure on $X$ is equivalent to the existence of a $\text{Pin}^+$-structure on $\Sigma$ that extends to the 2-handles. Since such 2-handles are attached to $\Sigma \times \partial D^2$ with relative odd framing, this corresponds to a $\text{Pin}^+$-structure on $\Sigma$ restricting to the one bounding the Möbius strip on every vanishing cycle. The conclusion follows from \cite{[D]}.

     In the $\text{Pin}^+$-case, $H^1(\Sigma;\Z_2)$ acts on the set of $\text{Pin}^+$-structures of the surface $\Sigma$ by
    \begin{equation}\label{act2}
        q^{+}_{\gamma}(x)=q^{+}(x)+ \gamma(x)
    \end{equation}
    for all $\gamma \in H^1(\Sigma;\Z_2)$ and $x \in H_1(\Sigma;\Z_4)$, see \cite{[D]}. In particular, $\text{Pin}^+$-structures on $\Sigma$ equivariantly correspond to quadratic enhancements and the action (\ref{act2}) is well defined on $H^1(X;\Z_2)\subset H^1(\Sigma;\Z_2)\cong H^1(X;\Z_2)\oplus K$, where $K \subset H^1(\Sigma;\Z_2)$ is the kernel of the map induced by the inclusion $\Sigma \subset X$.
\end{proof}

Let $f: X \to D^2$ be a Lefschetz fibration over the 2-disk with surface fiber $\Sigma$.  Let $c_1, \dots, c_n \subset \Sigma$ be its vanishing cycles. Note that one can find simple closed curves $e_1, \dots, e_g$ inducing a basis of $H_1(\Sigma;\Z_2)\cong \Z_2^g$. From now on, we will assume that $\Sigma$ has a $\text{Pin}^+$-structure, since this is a necessary condition for $X$ to be a $\text{Pin}^+$-manifold. In particular, this means that $\Sigma$ is not a closed non-orientable surface of odd Euler characteristics. Let
\begin{equation}\label{map}
    q^{+}_0: H_1(\Sigma;\Z_4) \to \Z_2
\end{equation}
to the quadratic enhancement associated to a fixed $\text{Pin}^+$-structure on $\Sigma$ and let $c_j=\sum_{i=1}^r c_{ji}[e_i].$ This gives the setting of Theorem \ref{Main2}, of which we now provide a proof.

\begin{proof}[Proof of Theorem \ref{Main2}]
    Lemma \ref{lemma2} implies that $X$ has a $\text{Pin}^+$-structure if and only if there is a map
    \[q^+: H_1(\Sigma,\Z_4) \rightarrow \Z_2\] 
    such that $q^+(x+y)=q^+(x)+q^+(y)+x \cdot y$ for all $x, y \in H_1(\Sigma; \Z_4)$ and $q^+([c_i])=1$ for all $i=1, \dots, n$. In particular, every such $q^+$ is necessarily of the form $q^{+}=q^{+}_0+l$ for some linear map $l=\sum_{i=1}^r \lambda_i [e^i]$, where $[e^i]\in H^1(\Sigma;\Z_2)$ denotes the dual element to $[e_i]\in H_1(\Sigma;\Z_2)$. In particular, this reduces the problem to finding $l= \sum_{i=1}^r \lambda_i e^i$ such that for every $j=1, \dots, n$ we have \[ q^{+}([c_j])=q^{+}_0([c_j])+l([c_j])=q_0^+([c_j])+\sum_{i=1}^r c_{ji} \lambda_i=1.\]
    The conclusion follows from Rouché-Capelli's Theorem \cite[Theorem 2.38]{[Linear]}.
\end{proof}
\begin{remark}
    The proof of Theorem \ref{theorem1} as well as the one of \cite[Theorem 1.1]{[Stipsicz]} rely on the fact that $H_1(\Sigma; \Z_2)$ is a vector space. This is why we can not adopt the same strategy for Theorem \ref{Main2}, since $H_1(\Sigma;\Z_4)$ is just a $\Z_4$-module. 
\end{remark}

\section{Some examples}\label{examples}
In this section, we show how to apply the results proven so far. In particular, starting from the handlebody decomposition of some small non-orientable 4-manifold $M$, we endow the union of handles up to index 2 with the structure of a Lefschetz fibration over the 2-disk, following the procedure explained in the proof of \cite[Theorem 1.1]{[MillerOzbagci]}. At this point, we apply Theorem \ref{theorem1} and Theorem \ref{Main2} to understand whether or not $M$ admits a $\text{Pin}^+$ or a $\text{Pin}^-$-structure, describing the action of $H^1(M;\Z_2)$ in terms of the action of a subgroup of $H^1(\Sigma;\Z_2)$ on the quadratic enhancements on the non-singular fiber $\Sigma$ of the Lefschetz fibration. We refer the reader to \cite[Section 1.5]{[Akbulutbook]} for the conventions and backgound on Kirby diagrams of non-orientable 4-manifolds, see also \cite{[Akbulut1]}, \cite[Section 2.1]{[BaisTorres]}, \cite{[CesardeSa]} and \cite{[MillerNaylor]}. The procedure we describe in this section applies to any closed non-orientable 4-manifold represented by means of a Kirby diagram.

\begin{example}[$\R \mathbb{P}^4$]
Figure \ref{1} represents a Kirby diagram for the the standard a handlebody decomposition of $\R \mathbb{P}^4$ with a single $k$-handle for each $k=0, \dots, 4$. In particular, the union of all handles up to index 2 corresponds to a tubular neighbourhood $\nu(\R \mathbb{P}^2) \cong D^2 \simtimes \R \mathbb{P}^2$ of $\R \mathbb{P}^2$ (see \cite[Section 0]{[Akbulut1]} and \cite[Section 2.5]{[BaisTorres]}) and $\R \mathbb{P}^4$ is the result of attaching a 3-handle and a 4-handle to this handlebody. Recall that 3-handles and 4-handles need not be drawn by \cite{[CesardeSa]} and \cite{[MillerNaylor]}, see \cite[Example 7]{[BaisTorres]}.
\begin{figure}
\minipage{0.4\textwidth}
  \includegraphics[width=80mm, trim=0 600 0 0, clip]{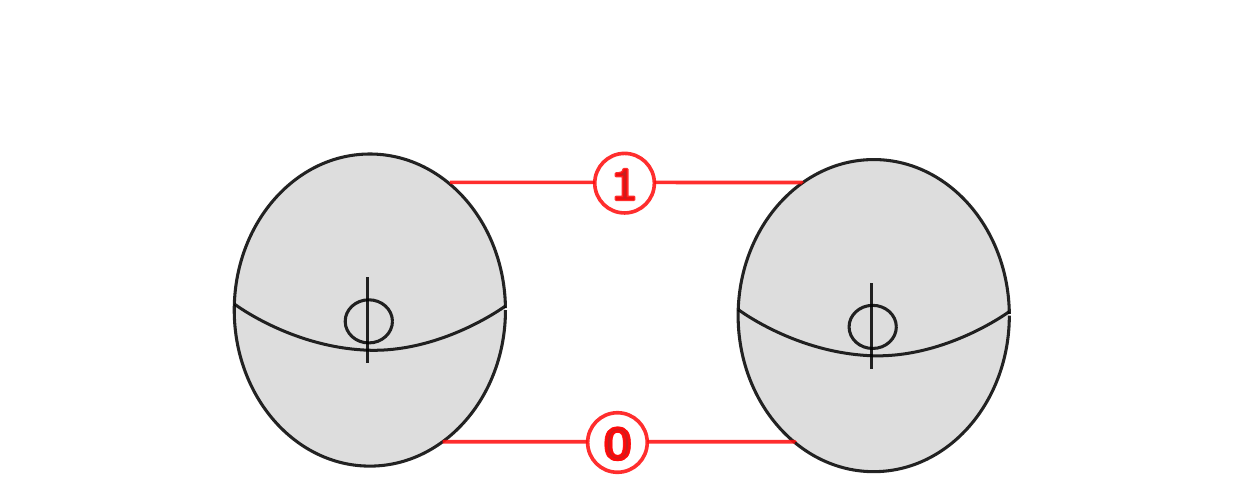} 
  \caption{A Kirby diagram of $\R \mathbb{P}^4$.}\label{1}
\endminipage\hfill
\minipage{0.6\textwidth}%
    \includegraphics[width=80mm, trim=0 500 0 50, clip]{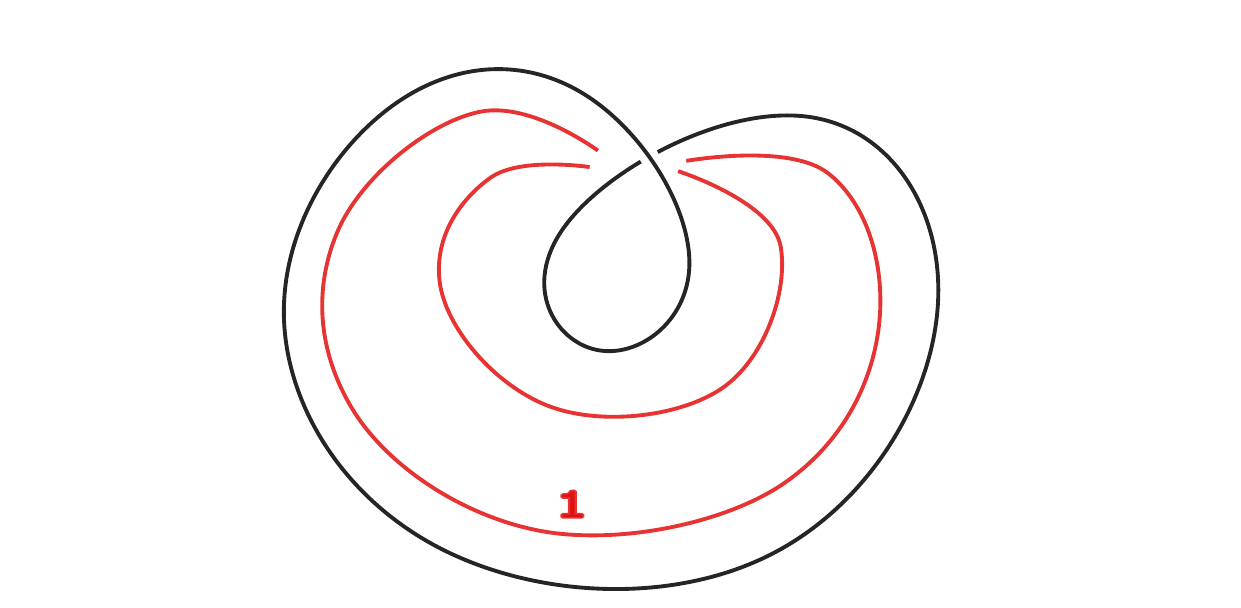} 
  \caption{The Lefschetz fibration associated to the 2-handlebody of $\R \mathbb{P}^4$.}\label{22}
\endminipage
\end{figure}

The union of handles up to index 2 can be seen as a Lefschetz fibration over the 2-disk with the Möbius band $\text{Mb}$ as fiber and a single vanishing cycle $c_1$, given by the projection of the $(1,0)$-framed 2-handle, see Figure \ref{22}.

A $\text{Pin}^+$-structure on $\text{Mb}$ extending to $\R \mathbb{P}^4$ is given by the quadratic enhancement

\begin{equation}
    q^{+}:H_1(\text{Mb};\Z_4) \rightarrow \Z_2
\end{equation}

defined by the condition $q^{+}([e_1])=1$, where $[e_1]$ is the generator of $H_1(\text{Mb};\Z_4) \cong \Z_4$ represented by the core of the band. Indeed, $[c_1]=2[e_1]$ and hence \[q^+([c_1])=q^{+}(2[e_1])= q^{+}([e_1])+q^{+}([e_1])+e_1^2=1 \in \Z_2.\] The other $\text{Pin}^+$-structure
\[q^{+}_x(y)=q^{+}(y)+x \cdot y \quad \text{ for all } y \in H_1(\text{Mb}; \Z_4)\]is obtained by acting on $q^{+}$ by the cohomology class dual to $[e_1]$, where $x \cdot y$ indicates the intersection number between $x$ and the $\Z_2$-reduction of $y$. Note that the condition $q^{+}_x([c_1])=1$ is still satisfied.

On the other hand, if $\R \mathbb{P}^4$ supported a $\text{Pin}^-$-structure, then there would be a quadratic enhancement
\[q^{-}:H_1(\text{Mb};\Z_2) \rightarrow \Z_4\]
satisfying $q^{-}([c_1])=2$. But this is impossible, since it would imply that
\[2=q^{-}([c_1])=q^{-}(2[e_1])=q^{-}(0)=0 \in \Z_4.\]
\end{example}

\begin{example}[$S^2 \simtimes \R \mathbb{P}^2$]
The non-orientable 4-manifold $S^2 \simtimes \R \mathbb{P}^2$ is obtained as a quotient of $S^2\times S^2$ under the orientation-reversing involution 
\begin{align*}
    \iota : S^2 \times S^2 \rightarrow S^2 \times S^2
    \\ (x,y) \mapsto (-x, \rho_{\pi}(y))
\end{align*}
where $\rho_{\theta}: S^2 \rightarrow S^2$ is the rotation of $S^2$ of angle $\theta$ about a fixed axis for any $\theta \in S^1$, see \cite[Chapter 12]{[Hillman]}. Moreover, it can also be seen as the result of a Gluck twist on the product $S^2 \times \R \mathbb{P}^2$ along a 2-sphere fiber, see \cite[Lemma 1]{[BaisTorres]}. In particular, there is a decomposition
\[ S^2 \simtimes \R \mathbb{P}^2= (S^2 \times D^2) \cup_{\phi} (S^2 \times \text{Mb})\] where the gluing map is \begin{align*}
\phi: S^1 \times S^1 \rightarrow S^2 \times S^2 \\ (x, \theta ) \mapsto (\rho_{\theta}(x), \theta).
\end{align*}

\begin{figure}
\minipage{0.4\textwidth}
  \includegraphics[width=80mm, trim=0 600 0 0, clip]{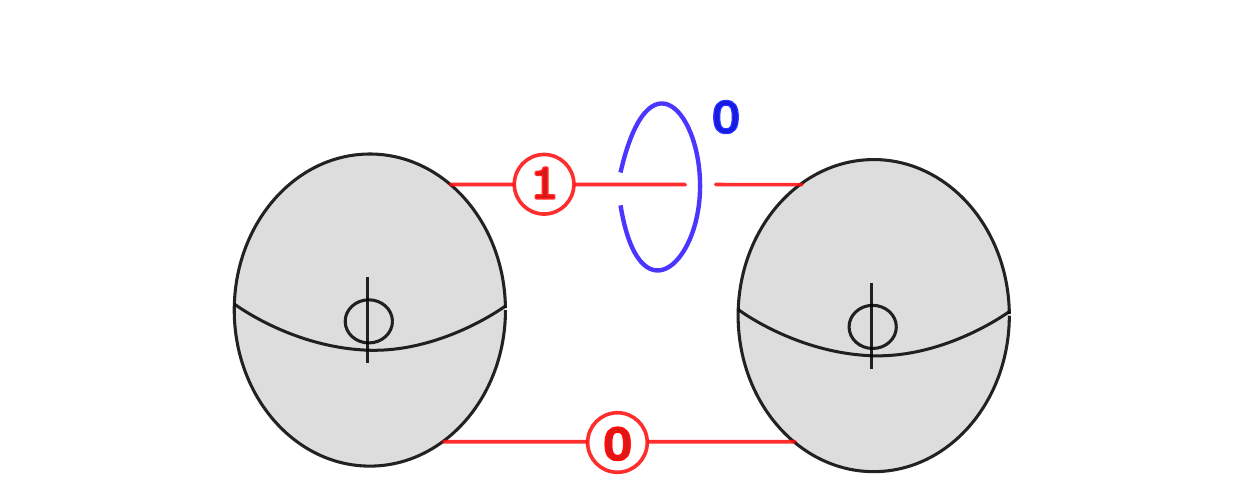} 
  \caption{A Kirby diagram of $S^2 \simtimes \R \mathbb{P}^2$.}\label{3}
\endminipage\hfill
\minipage{0.6\textwidth}%
   \includegraphics[width=80mm, trim=0 500 0 60, clip]{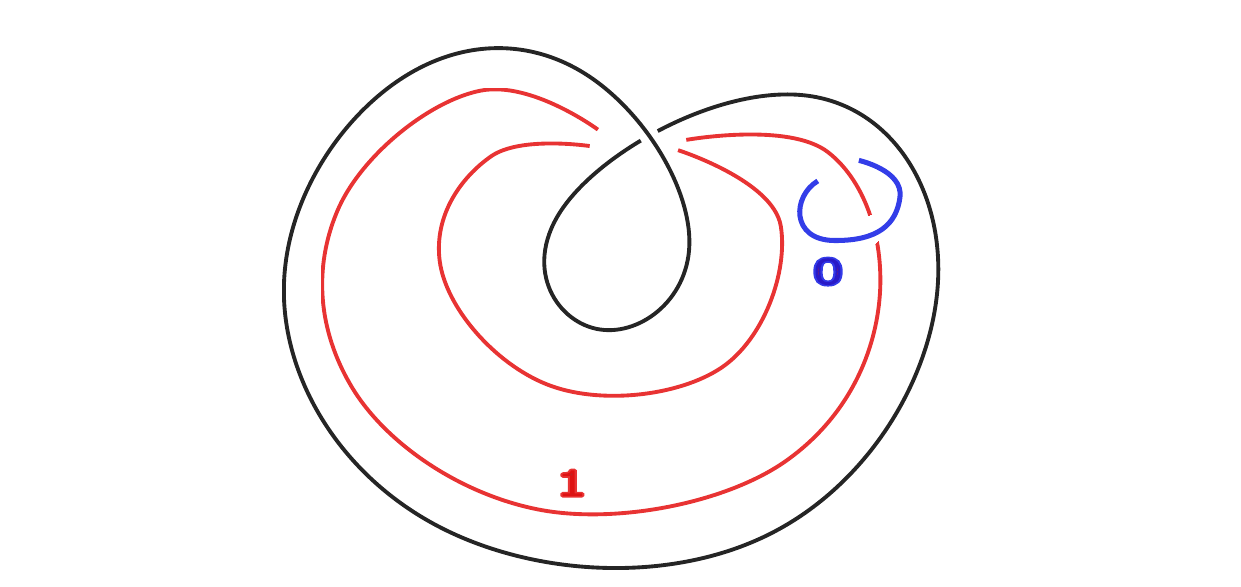} 
  \caption{}\label{4}
\endminipage
\end{figure}

A Kirby diagram of this manifold is given in Figure \ref{3}, see \cite[Section 0]{[Akbulut1]}, \cite[Figure 1]{[BaisTorres]}, \cite{[MillerNaylor]} and \cite[Section 4.1]{[Torres]}. It is obtained by considering $S^2  \simtimes \R \mathbb{P}^2$ as the double of $S^2 \simtimes D^2$. Figure \ref{4} represents the projection of the attaching circles of the 2-handles onto the page of the trivial Lefschetz fibration given by the union of handles up to index one. In order to endow the 2-handlebody of $S^2  \simtimes \R \mathbb{P}^2$ with the structure of a Lefschetz fibration over $D^2$, we need to stabilize the page and add vanishing cycles as in Figure \ref{5}, so that the attaching spheres of the 2-handles can be isotoped into different pages and have the right framing. This is Harer's trick and is explained in the non-orientable setting in \cite[Proof of Theorem 1.1]{[MillerOzbagci]}, see also \cite{[Harer]} and \cite[Theorem 2.1]{[EF]}. At the level of Kirby diagrams, this corresponds to adding canceling pairs of 1 and 2-handles.

\begin{figure}
\minipage{0.5\textwidth}
\includegraphics[width=70mm, trim=0 400 0 35, clip]{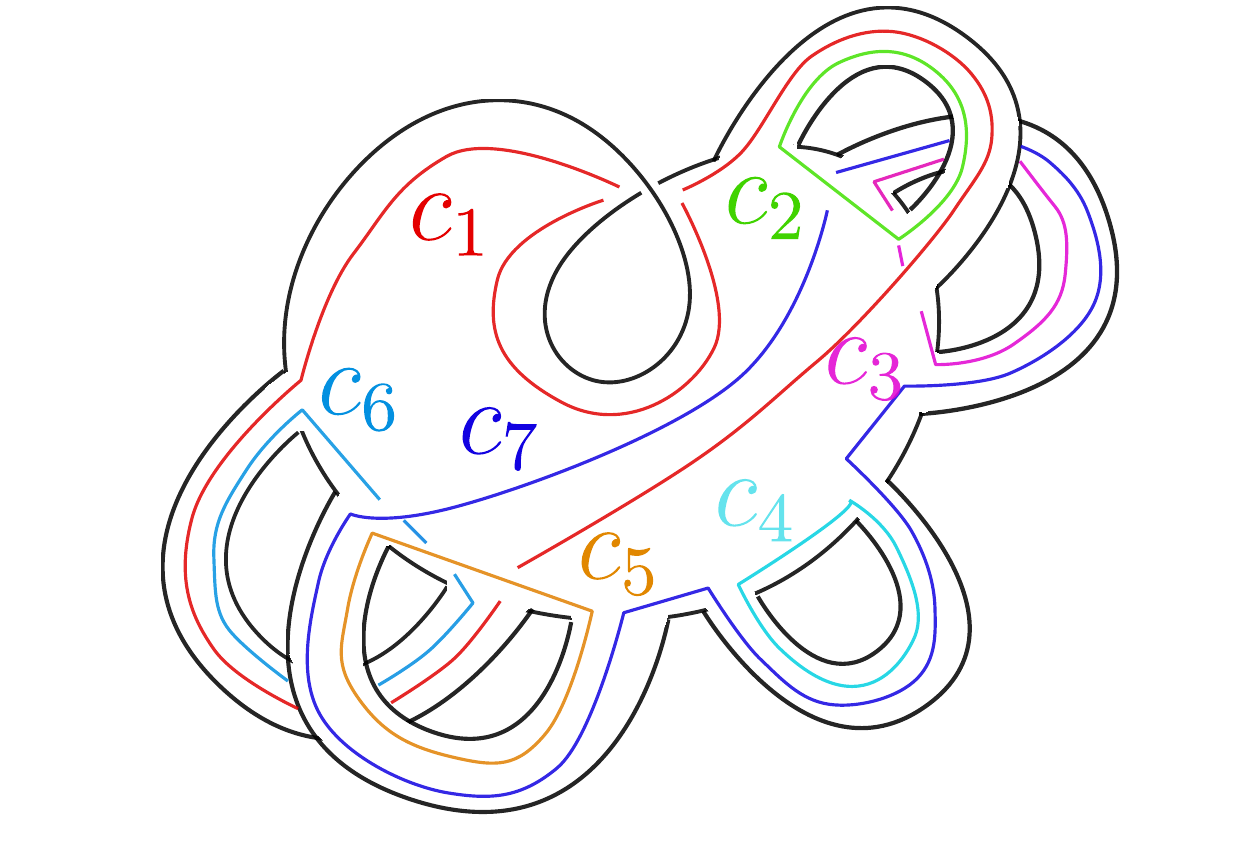} 
  \caption{The 2-handle-\ body of $S^2 \simtimes \R \mathbb{P}^2$ as a Lefschetz fibration over $D^2$.}\label{5}
\endminipage\hfill
\minipage{0.5\textwidth}%
     \includegraphics[width=70mm, trim=0 400 0 35, clip]{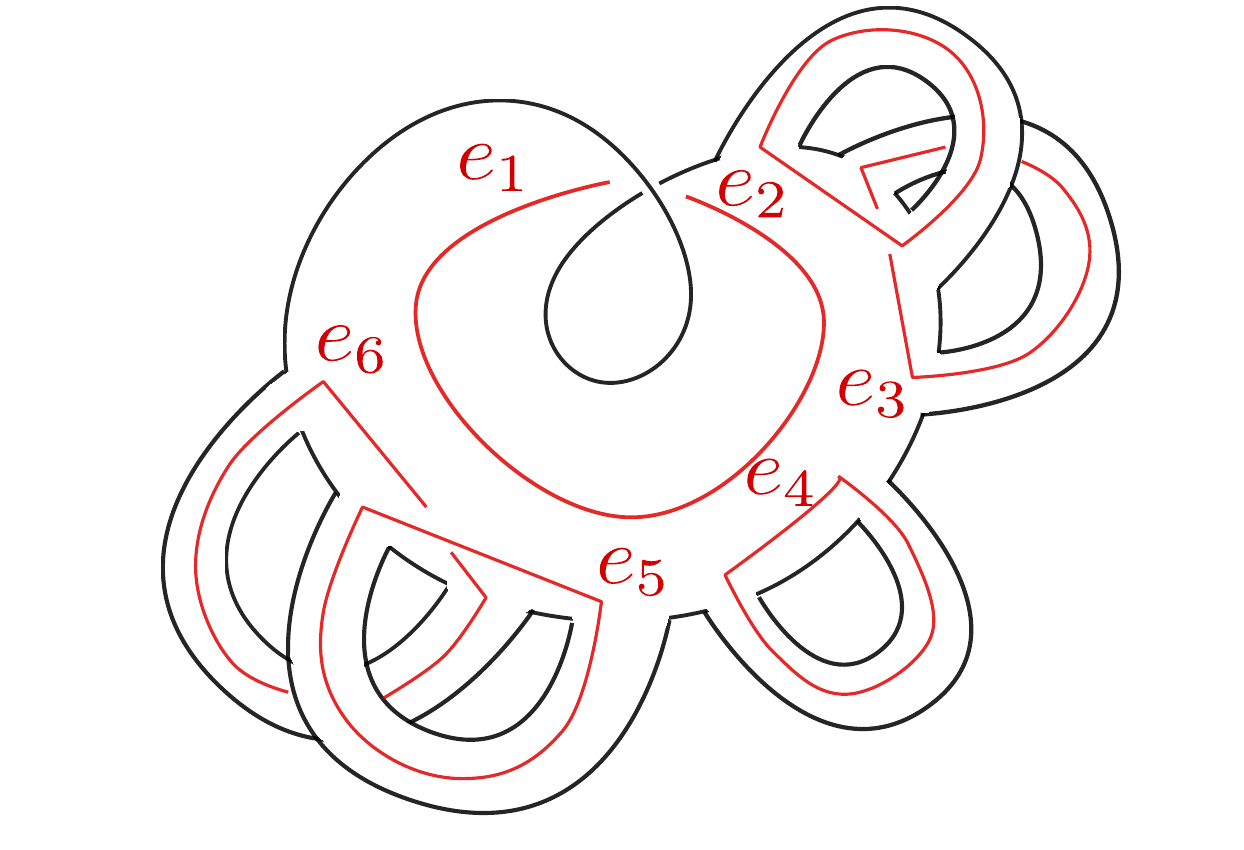} 
  \caption{Curves representing a basis of $H_1(\Sigma).$}\label{6}
\endminipage
\end{figure}

In Figure \ref{6} we fix a set of simple closed loops $e_1, \dots, e_6$ inducing a homology basis of the page $\Sigma$ of this new fibration. The quadratic enhancement
\[q^{+}:H_1(\Sigma;\Z_4) \rightarrow \Z_2\]
defined by the condition $q^{+}([e_i])=1$ for every $i=1, \dots, 6$ has the property that $q^{+}([c_i])=1$ for all the vanishing cycles $c_1, \dots, c_7$ and hence it determines a $\text{Pin}^+$-structure on $S^2 \simtimes \R \mathbb{P}^2$. The other $\text{Pin}^+$-structure can be found by acting on $q^{+}$ with $[e_1]\in H_1(\Sigma;\Z_2)$.

On the other hand, $S^2 \simtimes \R \mathbb{P}^2$ does not support a $\text{Pin}^-$-structure. Indeed, if this was the case, we could find a quadratic enhancement
\[q^{-}: H_1(\Sigma;\Z_2) \rightarrow \Z_4\]
such that $q^{-}([e_i])=2$ for all $i=1, \dots, 7$ and this would imply that
\[2=q^{-}([c_1])=q^{-}([e_2+e_6])=q^{-}([e_2])+q^{-}([e_6])=q^{-}([c_2])+q^{-}([c_6])=2+2=0\in \Z_4\] a contradiction.

\end{example}

\begin{example}[$S^2 \times \R \mathbb{P}^2$]

A Kirby diagram of the product 4-manifold $S^2 \times \R \mathbb{P}^2$ is given in Figure \ref{7}, see \cite[Section 0]{[Akbulut1]}, \cite[Figure 1]{[BaisTorres]}, \cite{[MillerNaylor]} and \cite[Section 4.1]{[Torres]}. It is obtained by considering $S^2 \times \R \mathbb{P}^2$ as the double of $D^2 \times \R \mathbb{P}^2$. Figure \ref{8} shows the projection of the attaching circles of the 2-handles onto the page of the trivial Lefschetz fibration given by the union of handles up to index one.

\begin{figure}
\minipage{0.4\textwidth}
   \includegraphics[width=80mm, trim=0 600 0 0, clip]{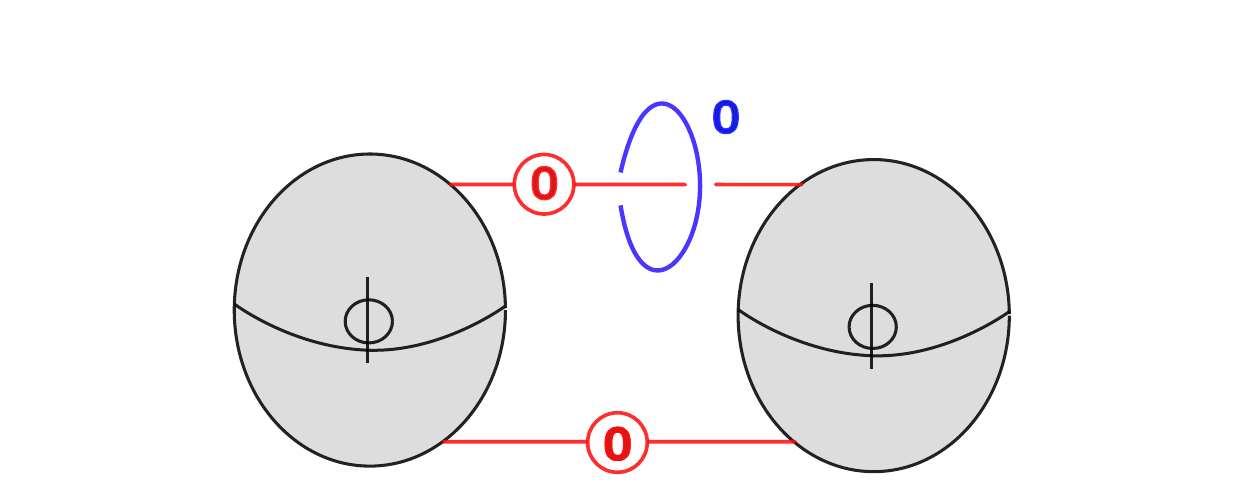} 
  \caption{A Kirby diagram of $S^2 \times \R \mathbb{P}^2$.}\label{7}
\endminipage\hfill
\minipage{0.6\textwidth}%
 \includegraphics[width=80mm, trim=0 500 0 60, clip]{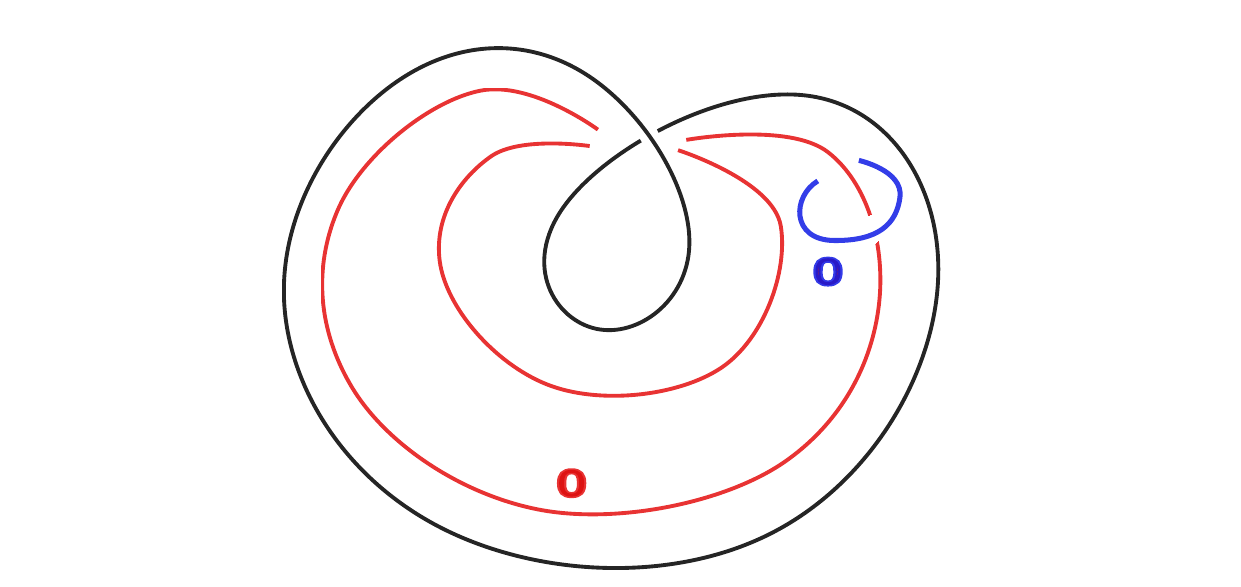} 
  \caption{}\label{8}
\endminipage
\end{figure}

  We use Harer'trick again (see \cite[Proof of Theorem 1.1]{[MillerOzbagci]}, \cite{[Harer]} and \cite[Theorem 2.1]{[EF]}) and show that the union of the handles up to index 2 in the fixed decomposition of $S^2 \times \R \mathbb{P}^2$ is given by the surface $\Sigma$ and by the vanishing cycles $c_1, \dots, c_8$ in Figure \ref{9}. Note that, with respect to the case of $S^2 \widetilde \times \R \mathbb{P}^2$, we need to stabilize the page one additional time in order to adjust the framing of the red colored 2-handle.

  \begin{figure}
\minipage{0.5\textwidth}
\includegraphics[width=70mm, trim=0 400 0 25, clip]{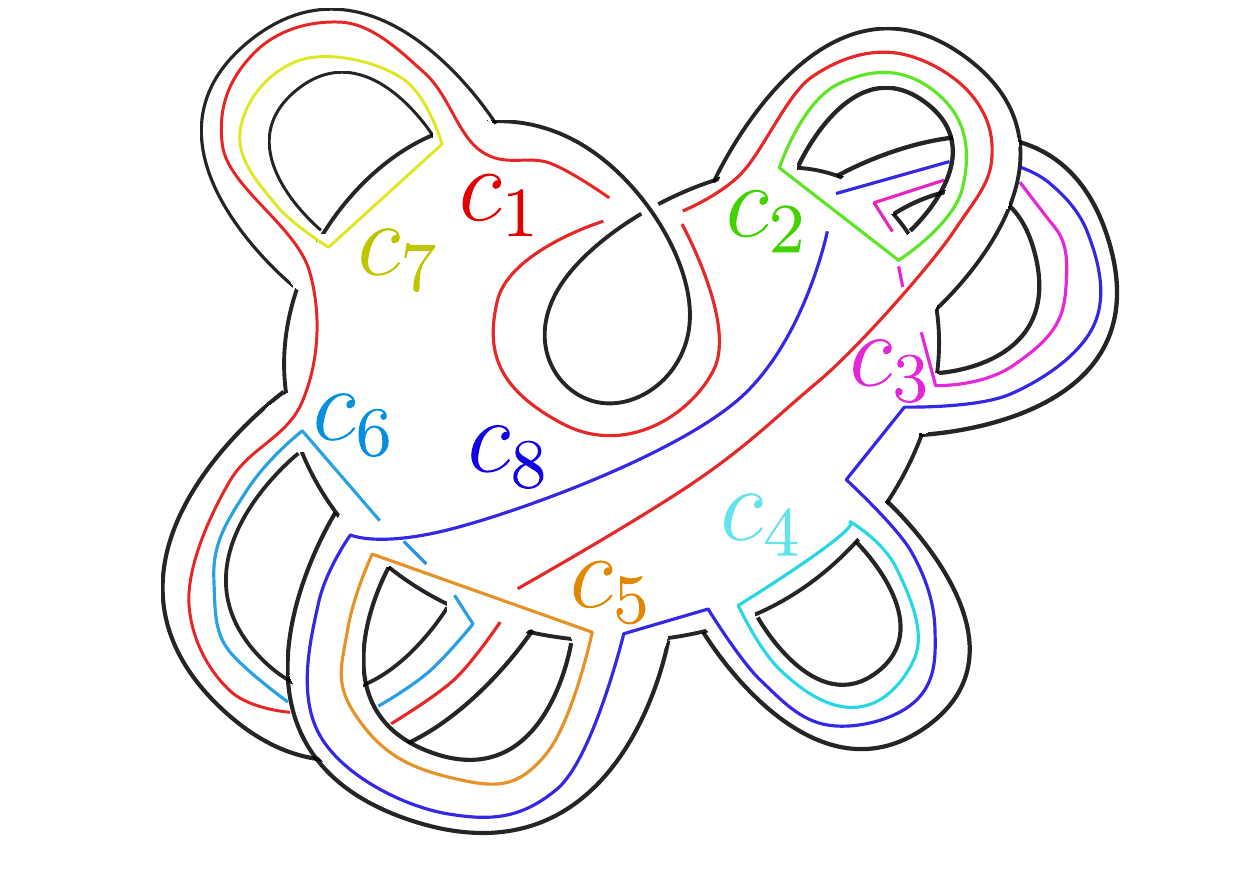} 
  \caption{The 2-handlebody of $S^2 \times \R \mathbb{P}^2$ as a Lefschetz fibration over $D^2$.}\label{9}
\endminipage\hfill
\minipage{0.5\textwidth}%
      \includegraphics[width=70mm, trim=0 400 0 25, clip]{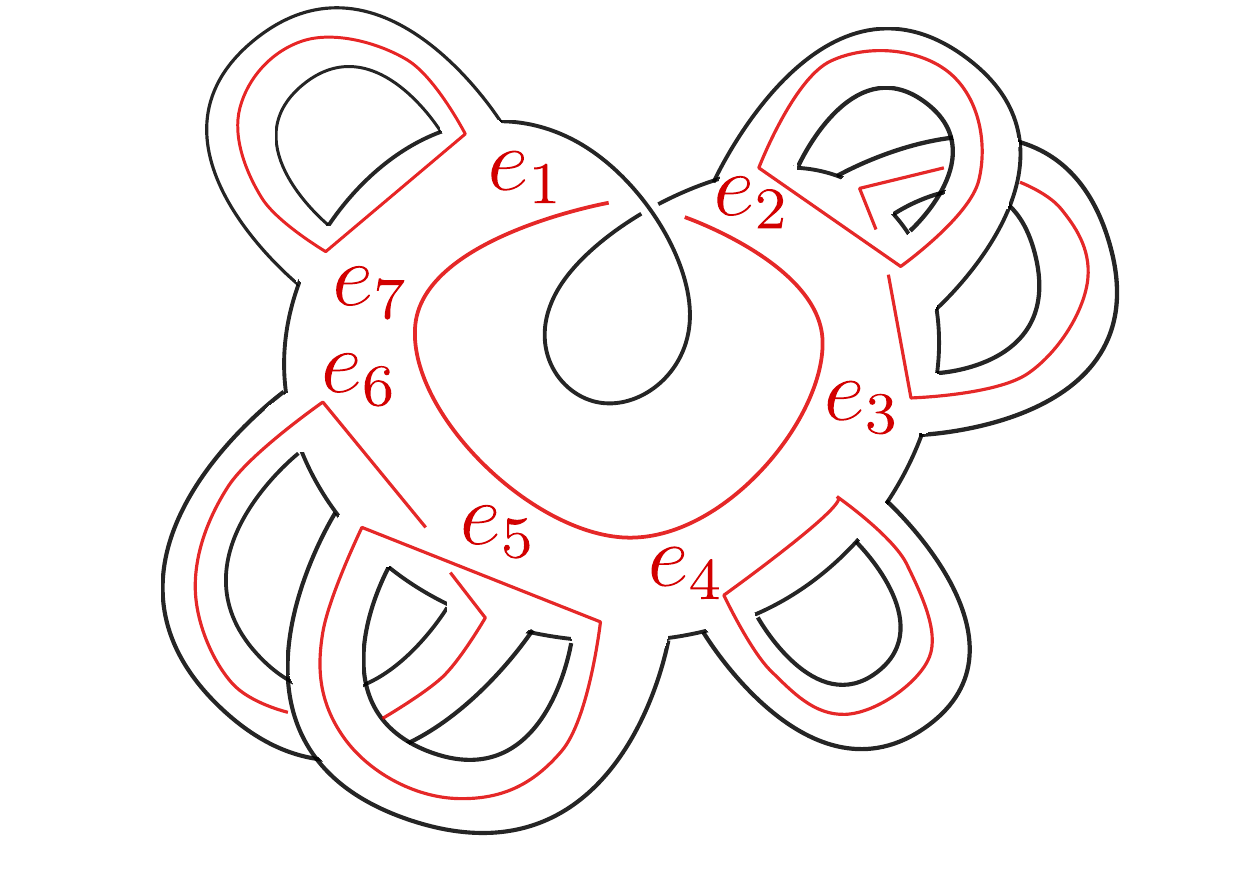}  
      \caption{Curves representing a basis of $H_1(\Sigma)$.}\label{10}
\endminipage
\end{figure}
 
  In Figure \ref{10} we again fix a set of loops $e_1, \dots, e_7 \subset \Sigma$ inducing a first homology basis of the page.

  The quadratic enhancement
\[q^{-}:H_1(\Sigma;\Z_2) \rightarrow \Z_4\]
defined by the condition $q^{-}([e_i])=2$ for every $i=1, \dots, 7$ has the property that $q^{-}([c_i])=2$ for all the vanishing cycles $c_1, \dots, c_8$ and hence it determines a $\text{Pin}^-$-structure on $S^2  \times \R \mathbb{P}^2$. The other $\text{Pin}^-$-structure can be found by acting on $q^{-}$ with $[e_1]\in H_1(\Sigma;\Z_2)$.

On the other hand, $S^2 \times \R \mathbb{P}^2$ can not support a $\text{Pin}^+$-structure. Indeed, if this was the case we could find a quadratic enhancement
\[q^{+}: H_1(\Sigma;\Z_4) \rightarrow \Z_2\]
such that $q^{+}([c_i])=1$ for all $i=1, \dots, 8$ and this would imply that
\begin{align*}
1=q^{+}([c_1])=q^{+}([2e_1+e_2+e_6+e_7])=\\  2q^{+}([e_1])+e_1^2+q^{+}([e_2])+q^{+}([e_6])+q^{+}([e_7])=\\ 1+q^{+}([c_2])+q^{+}([c_6])+q^{+}([c_7])=1+1+1+1=0 \in \Z_2     
\end{align*}
a contradiction.
\end{example}

\section{$\text{Pin}^{\pm}$-structures on Lefschetz fibrations over $S^2$}\label{sphere}

In the following, we provide necessary and sufficient conditions for the existence of a $\text{Pin}^{\pm}$-structure on a Lefschetz fibration over the 2-sphere. This is the non-orientable version of \cite[Theorem 1.3]{[Stipsicz]}.

\begin{theorem}
    Let $f: X \rightarrow S^2$ be a (possibly non-orientable) Lefschetz fibration with regular fiber $\Sigma$. Then $X$ supports a $\text{Pin}^-$-structure if and only if $X \setminus \nu(\Sigma)$ does and there exists a smoothly embedded surface $\sigma \subset X$ which is dual to $\Sigma$ in $H_2(X;\Z_2)$ and such that 
    \[[\sigma]^2+(w_1(\sigma)\cup w_1(\nu(\sigma))([\sigma])+w_1^2(\nu(\sigma))=0\in \Z_2.\]

    Analogously, $X$ supports a $\text{Pin}^+$-structure if and only if $X \setminus \nu(\Sigma)$ does and there exists a smoothly embedded surface $\sigma \subset X$ which is dual to $\Sigma$ in $H_2(X;\Z_2)$ and such that 
    \[\chi(\sigma)+[\sigma]^2+(w_1(\sigma)\cup w_1(\nu(\sigma)))([\sigma])=0\in \Z_2.\]
    
\end{theorem}
\begin{proof}
    The proof is a $\Z_2$-coefficients version of the one of \cite[Theorem 1.3]{[Stipsicz]} using the two formulas (\ref{-}) and (\ref{+}) in Section \ref{old}. We hence leave it to the interested reader as an exercise.
\end{proof}

\begin{remark}
    One can check whether $X \setminus \nu(\Sigma)$ supports a $\text{Pin}^{\pm}$-structure by applying Theorem \ref{theorem1} and Theorem \ref{Main2} respectively.
\end{remark}

\section {$\text{Pin}^+$ and $\text{Pin}^-$-structures on 3-manifolds}\label{3-maifolds}

The following well-known fact can be easily proven by adapting to the $\text{Pin}^-$ case the one of \cite[Theorem 1.4]{[Stipsicz]}, see \cite{[KirbyTaylor]}.

\begin{theorem}
    Any closed 3-manifold $M$ admits a $Pin^-$-structure.
\end{theorem}
\begin{proof}[Sketch of Proof]
    We reduce to the case in which $M$ is non-orientable, since if $M$ is orientable then it is automatically $\text{Spin}$ and hence also $\text{Pin}^+$ and $\text{Pin}^-$. $\R \mathbb{P}^2 \times S^1 \cs N$ for any 3-manifold $N$ is such an example. The idea is to write $M$ as the union of a non-orientable 3-dimensional handlebody $H$ of genus $g$ with $g$ 2-handles and a single 3-handle. Then \cite[Corollary 1.12]{[KirbyTaylor]} implies that $M$ admits a $\text{Pin}^-$-structure if and only if the closed surface $\partial H$ supports a $\text{Pin}^-$-structure which restricts to the one bounding the 2-disk on all the attaching spheres of the 2-handles and on all the belt spheres of the 1-handles. Such $\text{Pin}^-$-structure can always be constructed by means of a quadratic enhancement $q^{-}: H_1(\partial H;\Z_2) \rightarrow \Z_4$ which vanishes on all such curves. The proof follows the same lines of the one of \cite[Theorem 1.4]{[Stipsicz]} and we leave the details to the interested reader.
\end{proof}

On the other hand, it is not true that any closed 3-manifold admits a $\text{Pin}^+$-structure, since the condition $w_2(M)=0$ is not always satisfied in the non-orientable setting. However, we now show that it is still possible to check this condition by considering a handlebody decomposition of $M$.

\begin{theorem}
    Let $M=H \cup H'$ be a closed non-orientable 3-manifold obtained by adding $g$ 2-handles and a 3-handle to a genus-$g$ non-orientable handlebody $H$ and let $q_0^+: H_1(\partial H;\Z_4) \rightarrow \Z_2$ be a quadratic enhancement corresponding to a $\text{Pin}^+$ structure on $\partial H$. Let $a_1, \dots, a_g \in H_1(\partial H;\Z_4)$ and $b_1, \dots, b_g \in H_1(\partial H;\Z_4)$ be the homology classes of the attaching circles of the 2-handles and of the belt circles of the 1-handles of $H$ respectively and set \[q^{+}_0([a_j])=\alpha_j \quad \text{and} \quad q^{+}_0([b_j])=\beta_j\]
    for $j=1, \dots, g$. If for a fixed basis $e_1, \dots, e_{2g}$ of $H_1(\partial H;\Z_2)$ we have
    
    \[\tilde a_j=\sum_{i=1}^{2g} a_{j,i} e_i \quad \text{and} \quad  \tilde b_j=\sum_{i=1}^{2g} b_{j,i} e_i\]  for $j=1, \dots, g$, where $\tilde a_j$ and $\tilde b_j$ are the $\Z_2$-reduction of $a_j$ and $b_j$ respectively, then $M$ has a $\text{Pin}^+$-structure if and only if 

\[\text{rank}(C|A)=\text{rank} (C)\]
    where $C$ and $A$ are respectively the $2g \times 2g$ matrix and the column vector
    
    \begin{equation}
    C=
    \begin{pmatrix}
a_{1,1}  & \dots & a_{1,2g}\\
\vdots & \ddots &  \vdots\\
 a_{g,1} & \dots & a_{g,2g} \\
 b_{1,1} &  \dots & b_{1,2g}\\
\vdots & \ddots &  \vdots \\
 b_{2g,1} & \dots & b_{g,2g}
\end{pmatrix} 
, \qquad A=
    \begin{pmatrix}
\alpha_1\\
\vdots\\
\alpha_g\\ \beta_1\\
\vdots\\
\beta_g
\end{pmatrix} .
\end{equation} 
\end{theorem}

\begin{proof}
There is a $\text{Pin}^+$-structure on $M$ if and only if the non-orientable closed surface $\partial H$ has a $\text{Pin}^+$-structure which extends to the attaching circles of the 2-handles and to the belt spheres of the 1-handles in $H$. In particular, this is equivalent to check whether or not there exists a map 
\begin{equation}
    q^{+}: H_1(\partial H; \Z_4) \rightarrow \Z_2
\end{equation}
satisfying the condition $q^+(x+y)=q^+(x)+q^+(y)+x \cdot y$ for all $x,y \in H_1(\partial H;\Z_4)$ and which vanishes on $a_1, \dots, a_g$ and on $b_1, \dots, b_g$. 

If we fix a quadratic enhancement
\begin{equation}\label{H}
    q^{+}_0: H_1(\partial H; \Z_4) \rightarrow \Z_2
\end{equation}
corresponding to a $\text{Pin}^+$-structure on $\partial H$, then all the other $\text{Pin}^+$-structures on $\partial H$ will be of the form $q^{+}=q^{+}_0+l$ for some linear map $l \in \text{Hom}(H_1(\partial H;\Z_4), \Z_2)$. It follows that $M$ supports a $\text{Pin}^+$-structure if and only if there is $l=\sum_{i=1}^{2g} x_i e^i\in \text{Hom}(H_1(\partial H;\Z_4),\Z_2)$ such that \[\sum_{i=1}^{2g} a_{j,i}x_i=\alpha_j \quad \text{and} \quad \sum_{i=1}^{2g}b_{j,i}x_i=\beta_j\] in $\Z_2$ for $j=1, \dots, g$, where $e^1, \dots, e^{2g}$ is the dual basis to $e_1, \dots, e_{2g} \in H_1(\partial H; \Z_2)$. The conclusion follows from Rouché--Capelli's theorem.

\end{proof}
\section{Interpretation of $\text{Pin}^+$-structures à la Milnor}\label{Milnor}
Kirby and Taylor showed in \cite{[KirbyTaylor]} that the sets of $\text{Pin}^+$ and $\text{Pin}^-$-structures on a fixed vector bundle $\xi$ are in one to one correspondence with $\text{Spin}$-structures on $\xi \oplus 3 \cdot \text{det}(\xi)$ and $\xi \oplus \text{det}(\xi)$ respectively. In this section, we remark that for $\text{Pin}^+$-structures one can get another equivalent characterization, which recalls the following one of $\text{Spin}$-structures by Milnor.

\begin{theorem}[Milnor \cite{[Milnor]}]
    There is a canonical bijection between the set of $\text{Spin}$-structures on a real rank-$k$ vector bundle $\xi: E \rightarrow M$ and the set of homotopy classes of trivializations of $\xi|_{M^1}$ extending to trivializations of $\xi|_{M^2}$, where $M^i$ denotes the $i^{th}$ skeleton of $M$. Moreover, such correspondence is equivariant with respect to the action of $H^1(M;\Z_2)$.
\end{theorem}
In a similar fashion, we prove the following result. We remark that this follows essentially from the discussion in \cite{[KirbyTaylor]}.

\begin{theorem}
    Let $\xi: E \rightarrow M$ be a real rank-$k$ vector bundle. There is a canonical bijection between the set of $\text{Pin}^+$-structures on $\xi$ and the set of homotopy classes of $(k-1)$-tuples of everywhere linearly independent sections of $\xi|_{M^1}$ extending to $\xi|_{M^2}$, where $M^i$ denotes the $i^{th}$ skeleton of $M$. Moreover, such correspondence is equivariant with respect to the action of $H^2(M;\Z_2)$.
\end{theorem}
\begin{proof}
   The vanishing of $w_2(\xi)$ is a necessary and sufficient condition for the existence of both a $(k-1)$-tuple of linearly independent sections of $\xi|_{M^2}$ and a $\text{Pin}^+$-structure on $\xi$. We now suppose that $w_2(\xi)$ is trivial and we endow $\text{det}(\xi)$ with its canonical $\text{Pin}^+$-structure, see \cite[Addendum to 1.2]{[KirbyTaylor]}. Let $v_1, \dots, v_{k-1}$ be a $(k-1)$-tuple of linearly independent sections of $\xi|_{M^1}$ extending to $\xi|_{M^2}$. We are now going to associate to such a $(k-1)$-tuple a $\text{Pin}^+$-structure on $\xi|_{M^2}$. This will uniquely determine a $\text{Pin}^+$-structure on $\xi$, as in the case of $\text{Spin}$-structures. Using $v_1, \dots, v_{k-1}$ one can define a vector bundle isomorphism 
    \begin{equation}\label{eq}
    \xi|_{M^2} \cong \epsilon^{k-1} \oplus \text{det}(\xi)
    \end{equation}
    where $\epsilon^{k-1}$ denotes the trivial rank-$(k-1)$ real vector bundle. The associated $\text{Pin}^+$-structure on $\xi|_{M^2}$ is given by the pull-back via (\ref{eq}) of the $\text{Pin}^+$-structure on $\epsilon^{k-1} \oplus \text{det}(\xi)$ given by the direct sum of the trivial $\text{Spin}$-structure on $\epsilon^{k-1}$ with the fixed $\text{Pin}^+$-structure on $\text{det}(\xi)$. It is now easy to verify that this correspondence satisfies all the desired properties.
\end{proof}

\end{document}